\documentclass{amsart}

\usepackage[all]{xy}
 \usepackage[usenames,dvipsnames]{pstricks}
 \usepackage{epsfig}
 \usepackage{pst-grad} 
 \usepackage{pst-plot} 
\usepackage{ebezier}
\usepackage{amssymb,amsmath}

\newtheorem{thrm}{Theorem}

\newtheorem{rem}[thrm]{Remark}
\newtheorem{prop}[thrm]{Proposition}

\newcommand{\IM}{\operatorname{Im}}
\newcommand{\Ham}{\operatorname{Ham}}

\newcommand{\Crit}{\operatorname{Crit}}
\newcommand{\sing}{\operatorname{sing}}
\newcommand{\Image}{\operatorname{Im}}
\newcommand{\Fix}{\operatorname{Fix}}

\newcommand{\PD}{\operatorname{PD}}
\newcommand{\PSS}{\operatorname{PSS}}

\newcommand{\Morse}{\operatorname{Morse}}

\def \Dj{\mbox{\raise0.3ex\hbox{-}\kern-0.4em D}}

\def\l{\left}
\def\r{\right}

\def\a{\alpha}
\def\pa{\partial}
\def\b{\beta}
\def\A{\mathcal A}
\def\R{\mathbb R}

\begin{document}

\title{Comparison of spectral invariants in Lagrangian and Hamiltonian Floer theory}

\author{Jovana \Dj ureti\'c, Jelena Kati\'c and Darko Milinkovi\'c}

\address{Matemati\v{c}ki fakultet, Studentski trg 16, 11000
Belgrade, Serbia}

\email{jovanadj@matf.bg.ac.rs, jelenak@matf.bg.ac.rs, milinko@matf.bg.ac.rs}

\thanks{This work is
partially supported by Ministry of Education and Science
of Republic of Serbia Project \#ON174034.}

\begin{abstract} We compare spectral invariants in periodical orbits and Lagrangian Floer homology case,
for closed symplectic manifold $P$ and its closed Lagrangian submanifolds $L$, when $\omega|_{\pi_2(P,L)}=0$,
and $\mu|_{\pi_2(P,L)}=0$. We define a product $HF_*(H)\otimes HF_*(H:L)\to HF_*(H:L)$ and prove subadditivity of invariants with respect to this product.
\end{abstract}
\maketitle


2010 MS classification: Primary 53D40, Secondary 53D12, 57R58, 57R17


\section{Introduction}

 In~\cite{Oh1,Oh2} Oh defined spectral invariants for the case of cotangent bundle $P=T^*M$ (and
the canonical Liouville symplectic form), where the action functional
$$a_H(x):=\int_x\theta-\int_0^1H(x(t),t)dt$$ is well defined. Let $HF_*(H:O_M)$ denote Lagrangian Floer homology of the pair $(O_M,\phi_H^1(O_M))$, where $\phi_H^1$ is a time-one-map generated by a Hamiltonian $H$. Denote by
$HF_*^{\lambda}(H:O_M)$ the filtrated homology defined via filtrated Floer complex:
$$CF_*^{\lambda}(H:O_M):=\mathbb Z_2\langle\{x\in\Crit a_H\mid a_H(x)<\lambda\}\rangle.$$
These homology groups are well defined since the boundary map preserves the filtration:
$$\pa:CF_*^{\lambda}(H:O_M)\to CF_*^{\lambda}(H:O_M),$$
due to well defined action functional that decreases along its ``gradient flows".
For a singular homology class $\alpha$ define
$$\sigma(\alpha,H):=
\inf\{{\lambda}\in\R\mid F_H(\alpha)\in\IM(\imath_*^{\lambda})\}$$ where
$$\imath_*^{\lambda}:HF_*^{\lambda}(H:O_M)
\to HF_*(H:O_M)$$ is the homomorphism induced by inclusion and
$$F_H:H_*(M)\to HF_*(H:O_M)
$$ is an isomorphism between singular and Floer homology groups. The construction for spectral invariants is done in~\cite{Oh1} in the case of conormal bundle boundary condition, and in~\cite{Oh2} for cohomology classes. This construction is based on Viterbo's idea for generating functions
defined in the case of cotangent bundle (see~\cite{V}).

It turned out that Oh's invariants and the Viterbo's ones, are in fact the same, see~\cite{M1,M2}.

Oh proved in ~\cite{Oh1} that these invariants are independent both on the choice of almost
complex structure $J$ (which interferes in the definition of Floer homology) and, after a certain normalization,
on the choice of $H$ as far as $\phi^1_H(L_0)=L_1$. Using these invariants $\sigma(\alpha,L_1):=\sigma(\alpha,H)$, Oh derived the
non--degeneracy of Hofer's metric for Lagrangian submanifolds.

Lagrangian spectral invariants $\sigma$ were also used in~\cite{M3,M4} in the characterization of geodesics in
Hofer's metric for Lagrangian submanifolds of the cotangent bundle via quasi--autonomous Hamiltonians.

In~\cite{L}, Leclercq constructed spectral invariants for Lagrangian Floer theory without the assumption
$\omega=-d\theta$. He considered the case when $L$ is a closed submanifold of a compact (or convex in infinity)
symplectic manifold $P$ and
$$\omega|_{\pi_2(P,L)}=0,\quad \mu|_{\pi_2(P,L)}=0,$$ where $\mu$ stands for Maslov index.
He used a module structure of Floer ho\-mo\-lo\-gy over a Morse homology ring
and Albers' Piunikhin--Salamon--Schwarz (we will also use the abbreviation $\PSS$)
isomorphisms (see below or~\cite{A})
to prove that, after a certain normalization, the spectral invariants also do not depend on any included choices, but only on $L$ and
$L':=\phi^1_H(L)$.

Schwarz defined similar invariants in the case of Floer theory for contractible periodic orbits in~\cite{Sc2}.
If $(P,\omega)$ is a symplectic manifold with $\omega|_{\pi_2(P)}=0$, then the action functional is well defined
as:
$$\A_H(a):=\int_{D^2}\bar{a}^*\omega-\int_0^1H(a(t),t)dt,$$
where $\bar{a}:D^2\to P$ is any extension of $a$ to the unit disc.
For both $\omega|_{\pi_2(P)}=0$ and  $c_1|_{\pi_2(P)}=0$
Schwarz defined symplectic invariants as:
$$\rho(\a,H):=\inf\{\lambda\in\R\mid  \PSS(\a)\in\Image(\imath_*^{\lambda})\}.$$ Here $\a\in H^*_{\sing}(P)$ is a
nonzero cohomology class, $\PSS$ is a Piunikhin--Salamon--Schwarz isomorphism and $\imath_*^{\lambda}$
is a homomorphism induced by inclusion
$$\imath^{\lambda}:CF_*^{\lambda}(H)\to CF_*(H),
$$
where $CF_*(H)$
and $CF^{\lambda}_*(H)$ are (filtered) Floer chain complexes for Hamiltonian periodic orbits. For each nonzero
cohomology class $\a$, $\rho(\a,\cdot)$ is a section of the action
spectrum bundle
$$\begin{array}{l}\Sigma:=\bigcup_{\widetilde{\Phi}\in\widetilde{\Ham}(P)}\{\widetilde{\Phi}\}
\times\l\{a_H(x)\mid x\in\Fix^0\l(\Phi^1_H\r),\,\Phi^1_H\in\widetilde{\Phi}\r\}\\
\downarrow\\
\Ham(P,\omega)\end{array}$$ which is continuos with respect to Hofer's metric,
and which carries certain pro\-per\-ties (see~\cite{Sc2} for details and also~\cite{HZ}).

In his parer~\cite{A}, Albers constructed $\PSS$ morphisms for Lagrangian case and showed that, in certain dimensions, these
morphisms
are isomorphisms (see also~\cite{KM} for the case of cotangent bundle $P=T^*M$).
 These are the isomorphisms that Leclercq used in his already mentioned paper~\cite{L} to define the
Lagrangian invariants $\sigma$ (see Subsection~\ref{subsec_inv} below).
Albers considered the case of closed, monotone Lagrangian submanifold $L$ with minimal Maslov index
$\Sigma_L\ge 2$.
He also constructed morphisms
$\chi:HF_*(H:L)\rightarrow HF_*(H)$ and $\tau:HF_*(H)\to HF_*(H:L)$
where $HF_*(H)$ denotes Floer
homology for periodic Hamiltonian orbits and $HF_*(H:L)$ denotes Lagrangian Floer homology of the pair $(L,\phi_H^1(L))$. This construction is based on counting the numbers of ``chimneys" and it
was independently considered by
Abbondandolo and Schwarz in~\cite{AS}. Using these homomorphisms, Albers proved the commutativity
of certain diagrams (see~(\ref{diagA}) below).

In this paper, we consider the case of closed symplectic manifold $(P,\omega)$ and its closed smooth Lagrangian submanifold $L$ with topological assumptions
$$\omega|_{\pi_2(P,L)}=0,\quad \mu|_{\pi_2(P,L)}=0$$
and symplectic invariants $\rho$ for periodic orbits, and $\sigma$ for Lagrangian case (see Subsection~\ref{subsec_inv} below). 

We will use the homomorphisms constructed by ``chimneys" to compare these spectral invariants. Similar comparison was made by Monzner, Vichery and Zapolsky in different context (see~\cite{MVZ}). Further, we define the product $\circ$ using perturbed pseudoholomorphic curves that connect Hamiltonian periodic orbits and Hamiltonian paths with Lagrangian boundary conditions. This product was previously defined, for example by Hu and Lalonde~\cite{HL}, in the more general context of monotone Lagrangians. The main result of the paper is the following.

\begin{thrm}\label{thm:product} Let $P$ be a closed symplectic manifold and $L\subset P$ its closed Lagrangian submanifold such that $\omega|_{\pi_2(P,L)}=0$, $\mu|_{\pi_2(P,L)}=0$. Let $H_j:P\times [0,1]\to \R$ be three (time dependent) Hamiltonians, for $j=1,2,3$. Then there exists a product
$$\circ:HF_*(H_1)\otimes HF_*(H_2:L)\to HF_*(H_3:L)$$ which, in the case when $H_2=H_3$, turns the
Lagrangian Floer homology $HF_*(H_2:L)$  into a module over Floer homology for periodic orbits $HF_*(H_1)$. For $H_3=H_1\sharp H_2$, and $a\in HF_*(H_1)$, $b\in HF_*(H_2:L)$, it holds:
\begin{equation}\label{eq:ineq_circ}
\sigma(\PSS^{-1}(a\circ b),H_1\sharp H_2)\le\rho(\PSS^{-1}(a),H_1)+\sigma(\PSS^{-1}(b),H_2).
\end{equation}
\end{thrm}

The proof of Theorem~\ref{thm:product}  is given in Section~\ref{sec:proof_thm:product}.

In Section~\ref{sec:prel} we recall the construction on Floer homology and $\PSS-$type isomorphisms and their properties that we will use in the paper.

In Section~\ref{sec:inv} we construct spectral invariants in periodic orbits and Lagrangian case, and prove that Lagrangian spectral invariants do not depend on $H$ as long as $\phi^1_H(L)$ is fixed, up to a constant. Besides, we prove certain inequalities between these two types of invariants (Theorem~\ref{1.ineq} and Theorem~\ref{2.ineq}).

We would like to thank R\'emi Leclerq for pointing to us an error in the previous version of the paper.

We would like to thank the anonymous  referee for many valuable suggestions and corrections.

\section{Recalls and preliminaries}\label{sec:prel}

Throughout the paper we will assume that $(P,\omega)$ is a closed symplectic manifold and $L$ is its closed Lagrangian submanifold.

\subsection{Floer homology and $\PSS$ (iso)morphisms}\label{subsect:PSS}

Let us first briefly sketch the construction of Floer homology and $\PSS$ isomorphisms for periodical orbits. For a
smooth (generic) Hamiltonian $H:P\times S^1\to \R$,  Floer complex $CF_*(H)$ is defined as a vector space over
$\mathbb Z_2$ with the generators
$$\mathcal P(H):=\{a\in C^{\infty}(S^1,P)\mid \dot a(t)=X_H(a(t)), [a]=0\in\pi_1(P)\}$$
and it is graded by the Conley--Zehnder index (see~\cite{S}, for example). Floer differential is defined as
$$\delta(a):=\sum_bn(a,b)b,$$ where $n(a,b)$ is the number (modulo $2$) of the following set:
$$\mathcal M(a,b;H,J):=\left\{u:\R\times S^1\to P\left|
\begin{array}{ll}
\partial_s u+
J(\partial_t u-X_H(u))= 0\\
u(-\infty,t)=a(t)\\
u(+\infty,t)=b(t) \\
\end{array}\right.
\right\}$$ modulo $\R-$action $(\tau,u)\mapsto u(\cdot+\tau,\cdot)$.
Here $X_H$ is Hamiltonian vector field, i.e. $\omega(X_H,\cdot)=dH(\cdot)$.
Floer homology $HF_*(H)$ and Morse homology
$HM_*(P,f)$ for Morse function $f:P\to\R$ are  isomorphic. One way to prove this is via the $\PSS$
isomorphisms. Define two cut-off functions
$$\rho_R(s)=\begin{cases}1,&s\ge R+1,\\ 0, &s\le R\end{cases},\quad\tilde\rho_R(s):=\rho_R(-s).$$
Let $p$ be a critical point of $f$. Let
$g$ be a Riemannian metric such that the pair
$(f,g)$ is Morse--Smale.
Define
$$\PSS(p):=\sum_xn(p,x)x$$
and extend it on the chain level by linearity.
 Here $n(p,x)$ is a number (modulo $2$) of pairs $(\gamma,u):$
 $$\gamma:(-\infty,0]\to P,\quad u:[0,+\infty)\times S^1\to P$$ that satisfy
 $$\left\{
\begin{array}{ll}
\frac{d\gamma}{ds}=-\nabla f(\gamma(s))\\
\partial_s u+
J(\partial_t u-X_{\rho_RH}(u))= 0 \\
\gamma(-\infty)=p,\;
u(+\infty,t)=a(t)\\
\gamma(0)=u(0,\frac{1}{2}),\\
\end{array}
\right.$$
where $\nabla f$ denotes gradient of $f$ with respect to $g$.
This mapping commutes with the differentials, that is
$$\PSS\,\circ\,\partial_{\Morse}=\delta\circ \PSS,$$ where
$\partial_{\Morse}$ is Morse differential, so we have $\PSS:HM_*(P,f)\to HF_*(H)$ (we keep the same notation). It is proven in~\cite{PSS}
that, under our assumptions,
$\PSS$ is actually an isomorphism, and its inverse, $\PSS^{-1}$ is defined by counting the ``reverse" mixed object, i.e. the pairs
$(u,\gamma)$ that satisfy:
 $$\left\{
\begin{array}{ll}
u:(-\infty,0]\times S^1\to P, \quad \gamma:[0,+\infty)\to P\\
\partial_s u+
J(\partial_t u-X_{\tilde\rho_RH}(u))= 0 \\
\frac{d\gamma}{ds}=-\nabla f(\gamma(s))\\
u(-\infty,t)=a(t), \;
\gamma(+\infty)=p\\
u(0,\frac{1}{2})=\gamma(0)\\
\end{array}
\right.$$ (see Figure 1).

\begin{center}\includegraphics[width=12cm,height=2cm]{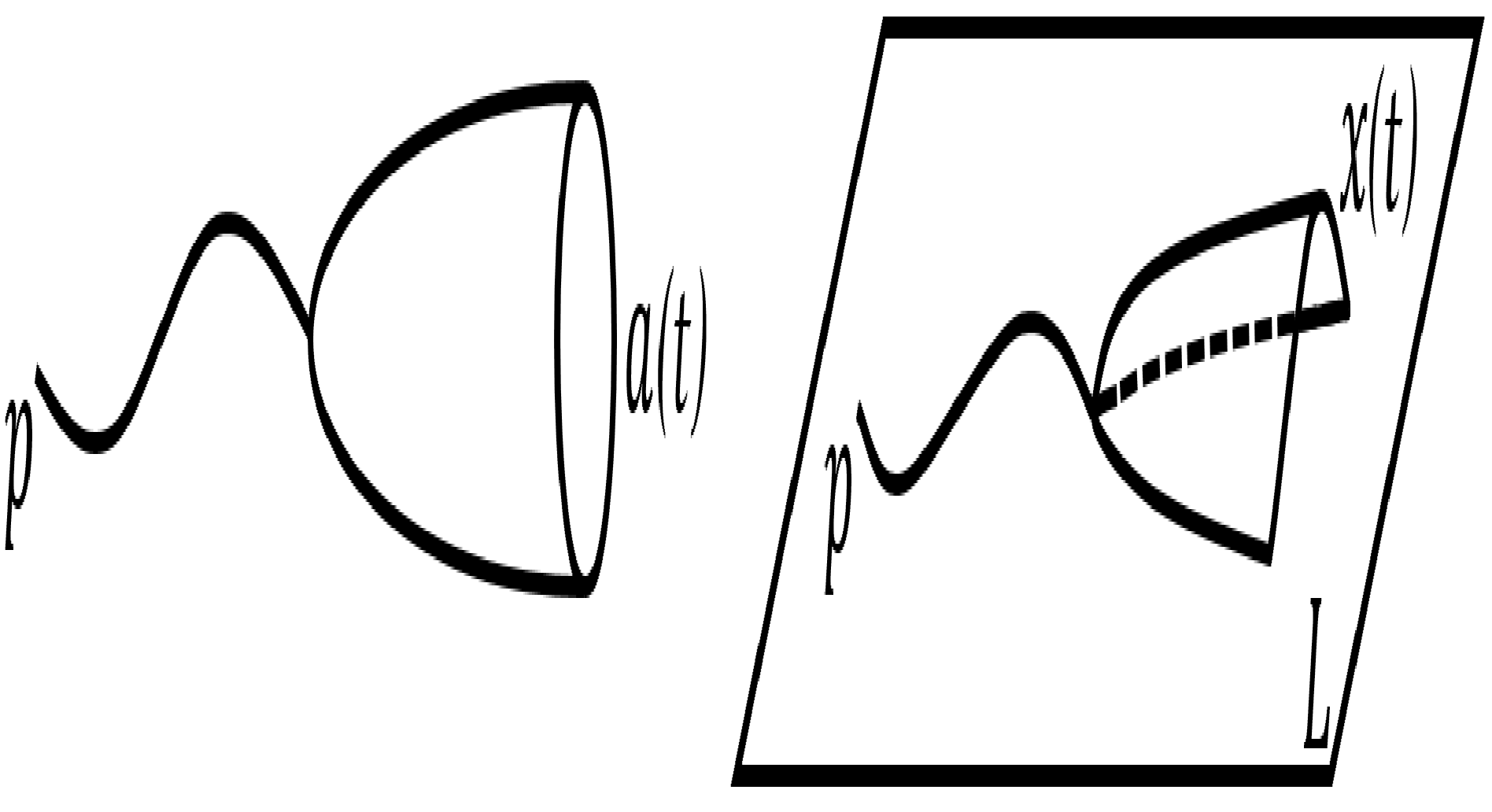}\\
\end{center}
\vspace{3mm}
\begin{center}{\bf Figure 1:}\, Mixed type objects that define PSS isomorphisms in the case of periodical orbits
and in Lagrangian intersections case
\end{center}

\medskip

Now let us recall the construction of Floer homology and $\PSS-$type morphisms in the Lagrangian case.
Suppose the $P$ and $L$ are as above and that $H$ is non degenerate, i.e. that $L\pitchfork \phi^1_H(L)$.
For our purposes, we will set
\begin{equation}\label{omega+mu}
\omega|_{\pi_2(P,L)}=0,\quad \mu|_{\pi_2(P,L)}=0
\end{equation} where $\mu$ is Maslov index. (Of course, Floer homology for Lagrangian intersections can be
defined in more general context.)
Let $H:P\times S^1\to\R$ be a smooth Hamiltonian function.
Floer complex $CF_*(H:L)$ is defined as a vector space over
$\mathbb Z_2$ with the generators
$$\mathcal P(H,L):=\{x\in C^{\infty}([0,1],P)\mid \dot x(t)=X_H(x(t)), x(0), x(1)\in L, [x]=0\in\pi_1(P,L)\}.$$
The grading is given by relative Maslov index,
which is well defined, since $ \mu|_{\pi_2(P,L)}=0$ (see, for example,~\cite{Oh} for details).
Floer differential is defined by counting the
pseudo--holomorphic tunnels, i.e.
$$\pa(x)=\sum_y n(x,y)y,$$ where $n(x,y)$ is the number (modulo $2$) of the set
$$\mathcal M(x,y;H,J):=\left\{u:\R\times[0,1]\to P\left|
\begin{array}{ll}
\partial_s u+
J(\partial_t u-X_H(u))= 0 \\
u(s, i)\in L, \; i\in\{0,1\}  \\
u(-\infty,t)=x(t) \\
u(+\infty,t)=y(t)  \\
\end{array}\right.
\right\}$$ modulo $\R-$action $(\tau,u)\mapsto u(\cdot+\tau,\cdot)$. As usual, we denote this quotient space by
$$\widehat{\mathcal M}(x,y;H,J):=\mathcal M(x,y;H,J)/\mathbb R.$$

The Albers' $\PSS-$type morphisms are well defined in more general cases than~(\ref{omega+mu}),
that is when $L$ is monotone and minimal Maslov number $N_L$ is at least $2$ (see~\cite{A}).
Let us recall this construction. For critical point $p$ of Morse function $f:L\to{\mathbb R}$, define
$$\Phi:CM_*(L,f)\to CF_*(H:L),\quad \Phi(p):=\sum_xn(p,x)x$$
 where $n(p,x)$ is a number (modulo $2$) of pairs $(\gamma,u)$, that satisfy
 \begin{equation}\label{moduli_PSS_A}\left\{
\begin{array}{ll}
\gamma:(-\infty,0]\to L,\quad u:[0,+\infty)\times [0,1]\to P\\
\frac{d\gamma}{ds}=-\nabla f(\gamma(s)) \\
\partial_s u+
J(\partial_t u-X_{\rho_RH}(u))= 0 \\
u(s,0),u(s,1),u(0,t)\in L \\
\gamma(-\infty)=p,\;
u(+\infty,t)=x(t) \\
\gamma(0)=u(0,\frac{1}{2})\\
\end{array}
\right.\end{equation}
(see Figure 1).
We denote the set of solutions of~(\ref{moduli_PSS_A}) by $\mathcal M_{p,x}^{f,H}$. The set
$\mathcal M_{p,x}^{f,H}$ is a $(m_{f}(p)-\mu(x))-$dimensional manifold where $m_f(p)$ is a Morse index of a critical point $p$
(note that this requires a particular choice of the reference of the Maslov index).

The map $\Phi$ turns out to be well defined in the homology level, and under our
assumption~(\ref{omega+mu}), an isomorphism between Morse and Floer  homologies in all dimensions.
We will denote this isomorphism of homology groups again by $\Phi$. Its inverse $\Psi$ is defined on the generators of Floer complex
as
$$\Psi:CF_*(H_1:L)\to CM_*(L,f), \quad\Psi(x):=\sum_pn(x,p)p$$
 where $n(x,p)$ is a number (modulo $2$) of pairs $(u,\gamma)$ that solve
 the equations:
 $$\left\{
\begin{array}{ll}
u:(-\infty,0]\times[0,1]\to P, \quad \gamma:[0,+\infty)\to L\\
\partial_s u+
J(\partial_t u-X_{\tilde\rho_RH}(u))= 0 \\
\frac{d\gamma}{ds}=-\nabla f(\gamma(s))\\
u(s,0),u(s,1),u(0,t)\in  L \\
u(-\infty,t)=x(t), \;
\gamma(+\infty)=p\\
u(0,\frac{1}{2})=\gamma(0).\\
\end{array}
\right.$$

The proofs of the above facts are usually based on the analysis of certain moduli spaces, especially in dimensions zero and one,
and their boundaries as well. The description of these boundaries and the proof of compactness  in zero dimension case use
Gromov compactness and gluing theorems. Bubbling is controlled due to topological assumptions~(\ref{omega+mu}).

For the sake of simplicity, we will denote these isomorphisms also by PSS, whenever there is no risk of confusion. More precisely
$$\PSS:=\Phi,\quad \PSS^{-1}=\Psi.$$

\section{Spectral invariants and their comparison}\label{sec:inv}

\subsection{Action functionals}\label{subsect:act_func}

In this subsection we will recall the constructions of two action functionals -- for contractible loops
and for contractible paths with the ends in Lagrangian submanifold.

In the case of periodic orbits, we will suppose that $\omega|_{\pi_2(P)}=0$, which is true if the second equality in~(\ref{omega+mu}) holds.
We define the
action functional on the space of smooth contractible loops
$$\Omega_0(P):=\{a\in C^{\infty}(S^1,P)\mid  [a]=0\in\pi_1(P)\}$$
in a standard way:
\begin{equation}\label{A_H}
\mathcal A_H(a):=-\iint_{D^2}\tilde{a}^*\omega-\int_{S^1}H(a(t),t)\,dt,
\end{equation}
where $\tilde{a}$ is any map from a disc with $\tilde{a}|_{S^1}=\gamma$. This map exists since $a$ is contractible
and the first integral in~(\ref{A_H}) does not depend on the choice of $\tilde{a}$ when the condition~(\ref{omega+mu}) is
fulfilled. One easily checks that the critical points of $\mathcal A_H$ are Hamiltonian periodic orbits.

Let us now define the action functional for Lagrangian case. Let $P$, $L$ and $H$ be as above and
suppose that
\begin{equation}\label{cond P}\omega|_{\pi_2(P,L)}=0.\end{equation}
The second condition in~(\ref{omega+mu}) does not have to be fulfilled in order to define the action functional, but it is necessary for Floer homology to be defined.
For the domain of the action functional $a_H$ we choose:
$$\Omega_0(P,L):=\{x\in C^{\infty}([0,1],P)\,\mid\,x(0), x(1)\in L,\,
[x]=0\in\pi_1(P,L)\}.$$ Set:
$$a_H(x,h):=-\iint_{D^2_+}h^*\omega
-\int_0^1H(x(t),t)\,dt,$$
where $h$ is any map from the upper half-disc $D^2_+$
to $P$ that restricts to $x$ on the upper half-circle.
Since $\omega|_{\pi_2(P,L)}=0$, the first integral does not depend on $h$, so we denote
$a_H(x):=a_H(x,h)$.

We compute the differential of $a_H$. For
 any variation $x_\varepsilon(t)$ of $x(t)$ with
$$x_\varepsilon(0),x_\varepsilon(1)\in L$$ let
$h_{\varepsilon}(s,t)$ be any smooth map from $D^2_+$ that satisfies
\begin{equation}\label{h_eps}\begin{aligned}&h_\varepsilon(0,t)\in L,\; \mbox{for}\; t\in[-1,1],\\
&h_{\varepsilon}(\cos(\pi\tau),\sin(\pi\tau))=x_{\varepsilon}(\tau)\; \mbox{for}\; \tau\in[0,1].
\end{aligned}\end{equation}
Denote by
$$\xi(t):=\frac{\pa}{\pa\varepsilon}\Big|_{\varepsilon=0}x_\varepsilon(t),\quad \zeta(s,t):=\frac{\pa}{\pa\varepsilon}\Big|_{\varepsilon=0}h_\varepsilon(s,t).$$
Using Cartan's and Stokes' formula, and the boundary conditions~(\ref{h_eps}),
one easily gets
$$\begin{aligned}
da_H(x)(\xi)=\frac{d}{d\varepsilon}\Big|_{\varepsilon=0}a_H(x_\varepsilon)=
-\iint_{D^2_+}\frac{d}{d\varepsilon}\Big|_{\varepsilon=0}
h_\varepsilon^*\omega-\int_0^1dH(\xi)dt=\\
-\iint_{D^2_+}d(i(\zeta)\omega)-\int_0^1dH(\xi)dt
=-\int_{\pa(D^2_+)}i(\zeta)\omega-\int_0^1dH(\xi)dt=\\
-\int_{S^1_+}\omega\l(\xi,\frac{dx}{dt}\r)-\int_0^1dH(\xi)dt
=-\int_0^1\omega\l(\xi,\frac{dx}{dt}-X_H\r)dt,
\end{aligned}$$
so the critical points of $a_H$ are Hamiltonian paths with ends in $L$.

\begin{rem}\rm If $y\in\Omega_0(P,L)\cap\Omega_0(P)$, i.e. $y(0)=y(1)\in L$ and $[y]=0\in\pi_1(P)$,
then $a_H(y)=\mathcal A_H(y)$.\qed
\end{rem}

\subsection{Invariants}\label{subsec_inv}

Now let us recall the definition of spectral invariants. We will start with periodic orbits case.
If~(\ref{omega+mu}) holds, we have well defined action functional $\mathcal A_H$.
Denote by
$$CF^{\lambda}_*(H):=\l\{
\sum c_aa\in CF_*(H)\mid c_a=0\;\mbox{for}\; \mathcal A_H(a)\ge\lambda\r\}.$$
Note that the Floer differential
$\delta$ preserves filtrations given by $\mathcal A_H$, and define
$$\delta_{\lambda}:=\delta|_{CF^{\lambda}_*(H)},\quad
HF^{\lambda}_*(H):=H_*(CF^{\lambda}_*(H),\delta_{\lambda}).$$
Denote by
$$\imath^{\lambda}_*:HF^{\lambda}_*(H)\to HF_*(H)$$ the homomorphism
induced by the inclusion map $\imath^{\lambda}$.
For $\a\in HM_*(P,f)$, define
\begin{equation}\label{eq:ro}
\rho(\a,H):=\inf \l\{\lambda\mid \,\PSS(\a)\in\Image\l(\imath^{\lambda}_*\r)\r\}.
\end{equation}
The above definition is also valid
in the case when $\a$ is a singular homological class, since Morse and singular homologies are isomorphic
(in the rest of the paper we will also sometimes identify Morse and singular homologies).

Now we consider the Lagrangian case. Suppose that $P$ and $L$
are closed and that they satisfy the condition~(\ref{omega+mu}). Suppose also that Hamiltonian
paths with the ends in $L$ belong to $\Omega_0$ (i.e. are zero in $\pi_1(P,L)$).
Since the action functional $a_H$ is well defined and since the differential $\pa$
preserves the filtration given by $a_H$, we can set
$$\begin{aligned}&CF^{\lambda}_*(H:L):=\l\{
\sum c_xx\in CF_*(H:L)\mid c_x=0\;\mbox{for}\; a_H(x)\ge\lambda\r\}\\
&\pa_{\lambda}:=\pa|_{CF^{\lambda}_*(H:L)}\\
&HF^{\lambda}_*(H:L):=H_*(CF^{\lambda}_*(H:L),\pa_{\lambda}).\end{aligned}$$
Denote by
$$\jmath^{\lambda}_*:HF^{\lambda}_*(H:L)\to HF_*(H:L)$$ the homomorphism
induced by the inclusion map $\jmath^{\lambda}$.
For given singular (or Morse) homological class $\a\in HM_*(L,f)$, define
$$
\sigma(\a,H):=\inf\l\{\lambda\mid \;\PSS(\a)\in\Image\l(\jmath^{\lambda}_*\r)\r\}.
$$

\begin{thrm}
If $\phi^1_H(L)=\phi^1_K(L)$, then
$$\sigma(\a,H)-\sigma(\a,K)=C=C(H,K).$$
\end{thrm}

\noindent{\it Proof:} Denote by $L_1=\phi^1_H(L)=\phi^1_K(L)$. Let $x, y\in L\cap L_1$. Let $c(s)$ be any smooth path in $L$ 
connecting $x$ and $y$ (s.t. $c(0)=x$, $c(1)=y$). Denote by
$$\begin{aligned}
&\gamma_{H,x}(t):=\phi^t_H(\phi^{-1}(x)),\quad \gamma_{K,x}(t):=\phi^t_K(\phi^{-1}(x)),\\
&\gamma_H^s(t):=\phi^t_H(\phi^{-1}(c(s))),\quad \gamma_K^s(t):=\phi^t_K(\phi^{-1}(c(s))),\\
&u_H(s,t):=\gamma_H^s(t),\quad u_K(s,t):=\gamma_K^s(t).
\end{aligned}$$
Define $f(s)$ to be 
$$f(s):=\iint_{D_s}\omega-\l(\int_{\gamma_H^s}Hdt-\int_{\gamma_K^s}Kdt\r),$$
where $D_s$ is the surface consisting of the union of $u_H(\tau,t)$, $u_K(\tau,t)$ $(\tau,t)\in[0,s]\times[0,1]$, and the two half discs with the boundaries on $\gamma_H$ and $L$ (respectively $\gamma_K$ and $L$), which exists due to the assumption $[\gamma_{H,x}]=[\gamma_{K,x}]=0\in\pi_1(P,L)$. Now, using Stokes' and Cartan's formula, as well as the condition $\pi_2(P,L)=0$, one easily derives $f'(s)=0$. Obviously, $f(0)=a_H(\gamma_{H,x})-a_K(\gamma_{K,x})$. Since $\pi_2(P,L)=0$, it holds $f(1)=a_H(\gamma_{H,y})-a_K(\gamma_{K,y})$. This means that $f\equiv const$, i. e. 
$$a_H(\gamma_{H,x})-a_K(\gamma_{K,x})=a_H(\gamma_{H,y})-a_K(\gamma_{K,y}).$$
We now proceed as in~\cite{Oh1}, namely, we switch to the geometric, instead of dynamic version of Floer homology. More precisely, there is a transformation between $L\cap L_1$ and Hamiltonian paths with ends on $L$, as well as perturbed holomorphic discs with boundary on $L$ on one side, and {\it holomorphic} discs but with the conditions $u(s,0)\in L$, $u(s,1)\in L_1$, so we have an one-to-one correspondence
between the generators and the boundary operators in two versions of Floer homology. Therefore, the set of elements of Floer homology that participate in the definitions on invariants are the same in two versions, so the claim follows. 
\qed

\subsection{Chimneys and comparison between invariants}

Let us first recall the Albers' construction of a homomorphism between Floer homology for periodical Hamiltonian orbits
and Lagrangian Floer homology. We assume that Hamiltonian
$H:P\times[0,1]\to\R$ is admissible in the sense of~\cite{A}, meaning that there are no constant contractible periodic orbits.
For
$$\Sigma:=\R\times [0,1]/\sim,\quad\mbox{where}\; (s,0)\sim (s,1)\;\mbox{for}\;s\le 0,$$
$a\in CF_*(H)$ and $x\in CF_*(H:L)$ define the manifold
of chimneys as:
$$\mathcal M(a,x):=\l\{u:\Sigma\to P\l|\begin{array}{l}\pa_su+J(\pa_tu-X_H\circ u)=0\\
u(s,0),u(s,1)\in L,\;\mbox{for}\;s\ge 0\\
u(-\infty,t)=a(t), \,u(+\infty,t)=x(t)
\end{array}\r.\r\}$$
(see Figure 2).
For $a\in CF_*(H)$, define
$$\tau(a):=\sum n(a,x)\,x$$ where
$n(a,x)$ stands for the number $\pmod 2$ of zero-dimensional component of $\mathcal M(a,x)$. This homomorphism
descends to $HF_*(H)$, namely, since $\tau\circ\delta=\pa\circ\tau$, it is well defined as a map:
$$\tau:HF_*(H)\to HF_*(H:L).$$ The following diagram commutes:
\begin{equation}\label{diagA}
\xymatrix{
  HF_*(H)\ar[d]_{\PSS^{-1}}\ar[r]^-{\tau} &HF_*(H:L)
  \ar[d]_{\PSS^{-1}}\\
H_*(P)\ar[r]^{\imath_!} &H_*(L)}
\end{equation}
where
 $\imath_!=\PD^{-1}\circ \,\imath_*\circ\PD$ is a homomorphism defined by Poincar\'e duality map and
 the inclusion and $H_*(P)$ and $H_*(L)$ are singular or Morse homologies
(see~\cite{A} for the details).

\begin{center}\includegraphics[width=7cm,height=3.5cm]{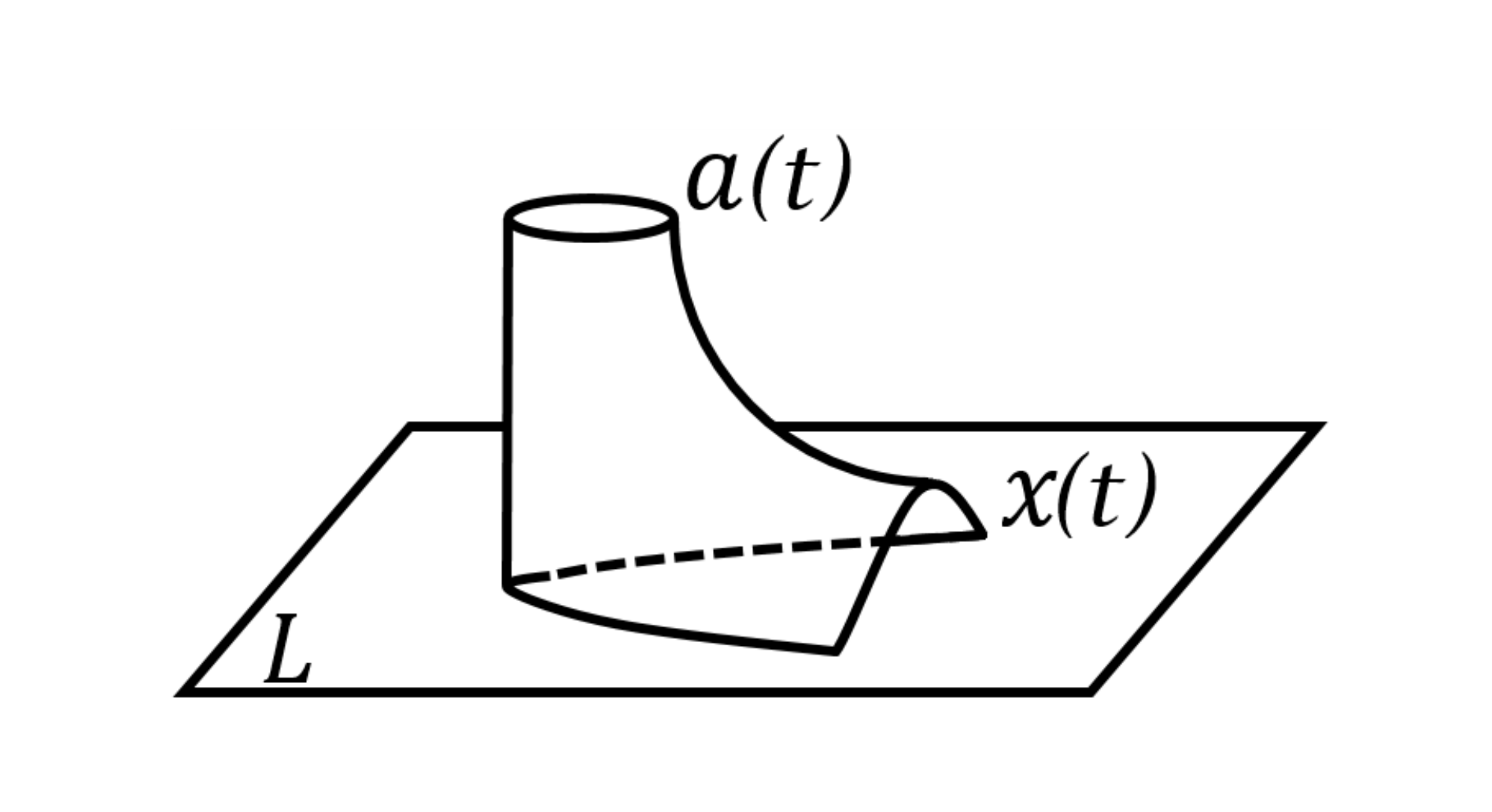}\\
\end{center}
\begin{center}{\bf Figure 2:}\, A "chimney" that defines the homomorphism $\tau$
\end{center}

\medskip

\begin{prop} Homomorphism $\tau$ induces a homomorphism $\tau_{\lambda}$ on filtered homology groups:
$$\tau_{\lambda}:HF^{\lambda}_*(H)\to HF^{\lambda}_*(H:L).$$
\end{prop}

\noindent{\it Proof:} Let $a\in CF^{\lambda}_*(H)$ and $x\in CF^{\lambda}_*(H:L)$
such that there exists $u\in \mathcal M(a,x)$. Denote by $y(t):=u(0,t)$. Since $a\in \Omega_0(P)$ and there exists
$u$ connecting $y$ and $a$, we have $y\in\Omega_0(P)$. Similarly, from $x\in \Omega_0(P,L)$ and the existence of
$u$ connecting $y$ and $x$, we conclude $y\in\Omega_0(P,L)$.
Since $y\in\Omega_0(P,L)\cap\Omega_0(P)$,
it holds $a_H(y)=\mathcal A_H(y)$, so we have:
\begin{equation}\label{tau_descends}
\begin{aligned}&a_H(x)-\mathcal A_H(a)=a_H(x)-a_H(y)+\mathcal A_H(y)-\mathcal A_H(a)=\\
&\int_0^{+\infty}\frac{d}{ds}a_H(u(s,\cdot))ds+
\int_{-\infty}^{0}\frac{d}{ds}\mathcal A_H(u(s,\cdot))ds=\\
&-\int_{-\infty}^{\infty}\int_0^1\omega\l(\pa_su,\pa_tu-X_H\circ u\r)\,dtds=\\
&-\int_{-\infty}^{+\infty}\int_0^1\l\|\pa_su\r\|^2_{g_J}\,dtds\le 0.
\end{aligned}\end{equation}
Hence, if $a\in CF^{\lambda}_*(H)$, then $\tau(a)\in CF^{\lambda}_*(H:L)$.
So, we can define
$$\tau_{\lambda}:=\tau|_{CF^{\lambda}_*(H)}:CF^{\lambda}_*(H)\to CF^{\lambda}_*(H:L).$$
Let $a$ be in $\Image(\delta_{\lambda})$, i.e. $a=\delta b$, with $b\in CF_*(H)$
and $\mathcal A_H(b)<\lambda$. We know that $\mathcal A_H(a)<\lambda$ also, since the action functional
decreases along perturbed holomorphic strips that define the differential $\delta$.
From~(\ref{tau_descends}) we have that $a_H(\tau(b))<\lambda$, so
$$\tau_{\lambda}(\delta_\lambda b)=\tau_{\lambda}(a)=\tau(a)=\tau(\delta b)=\pa\tau(b)
=\pa_{\lambda}\tau_{\lambda}(b),$$
since $\tau$ commutes with the differentials. This implies that $\tau_{\lambda}$
descends to the homology level.\qed

\begin{rem}\rm It is obvious that the diagram
\begin{equation}\label{diag i}\xymatrix{
  HF_*^{\lambda}(H)\ar[d]_{\imath_*^{\lambda}}\ar[r]^-{\tau_{\lambda}} &HF_*^{\lambda}(H:L)
  \ar[d]_{\jmath_*^{\lambda}}\\
HF_*(H)\ar[r]^-{\tau} &HF_*(H:L)}\end{equation} commutes.\qed
 \end{rem}

 \begin{thrm}\label{1.ineq} If $\a\in H_*(P)$ is a singular (or Morse) homological class, then $\rho(\alpha,H)\ge\sigma(\PSS(\imath_!(\a)),H)$.
 \end{thrm}
\noindent{\it Proof:} Consider the following commutative diagrams:
\begin{equation}\label{2diags}
\xymatrix{ HF_*^{\lambda}(H)\ar[d]_{\imath_*^{\lambda}}\ar[r]^-{\tau_{\lambda}} &HF_*^{\lambda}(H:L)
\ar[d]_{\jmath_*^{\lambda}}\\
HF_*(H)\ar[d]_{\PSS^{-1}}\ar[r]^-{\tau} &HF_*(H:L)
\ar[d]_{\PSS^{-1}}\\
H_*(P)\ar[r]^{\imath_!} &H_*(L)}
\end{equation}
The upper diagram is~(\ref{diag i}) and the lower is Albers'~(\ref{diagA}). For given $\a\in H_*(P)$ and
$\b\in H_*(L)$, let us define the sets:
$$\begin{aligned}
&A_H(\a):=\{\lambda\mid \,\PSS(\a)\in\Image(\imath^{\lambda}_*)\}\\
&A_{H:L}(\b):=\{\lambda\mid \;\PSS(\b)\in\Image(\jmath^{\lambda}_*)\}.
\end{aligned}
$$
Let $\lambda\in A_H(\a)$. There exists $a\in HF_*^{\lambda}(H)$ such that
$\PSS^{-1}(\imath^{\lambda}_*(a))=\a$. Since both diagrams~(\ref{2diags})
commute, this implies that
$\PSS^{-1}(\jmath^{\lambda}_*(\tau_{\lambda}(a)))=\imath_!(\a)$, so
$\PSS(\imath_!(\a))\in\Image(\jmath^{\lambda}_*)$. This means that
$\lambda\in A_{H:L}(\imath_!(\a))$ i.e.
$$A_H(\a)\subset A_{H:L}(\imath_!(\a)).$$
Since
$$\rho(\a,H)=\inf A_H(\a),\quad \sigma(\b,H)=\inf A_{H:L}(\b)$$ the claim follows.\qed

 In the same way, considering Albers' commutative diagram
 $$\xymatrix{
  HF_*(H:L)\ar[d]_{\PSS^{-1}}\ar[r]^{\chi} &HF_*(H)\ar[d]_{\PSS^{-1}}\\
H_*(L)\ar[r]^{\imath_*} &H_*(P)}$$ where $\chi$ is also defined using chimneys, but in opposite direction
(see~\cite{A}), one can prove the following

 \begin{thrm}\label{2.ineq} If $\b\in H_*(L)$ is a singular (or Morse) homological class, then
 $\rho(\imath_*(\b),H)\le\sigma(\b,H)$.
 \qed\end{thrm}

 \vspace{1cm}

\section{Proof of Theorem~\ref{thm:product}}\label{sec:proof_thm:product}

 The product
 $$\circ:HF_*(H_1)\otimes HF_*(H_2:L)\to HF_*(H_3:L)$$ is define by counting of a sort of pair-of-pants objects. More precisely, consider the disjoint union
 $$\mathbb R\times[-1,0]\sqcup\mathbb R\times[0,1]$$ and identify $(s,0^-)$ with $(s,0^+)$ for all $s\ge 0$ as well as $(s,0^+)$ with $(s,1)$, for $s\le 0$ (see figure below). Denote the obtained Riemannian surface with boundary by $\Sigma$. Denote by $\Sigma_1$, $\Sigma_2$, $\Sigma_3$ the two "incoming" and one "outgoing" ends, such that
 $$\begin{aligned}
 &\Sigma_1\approx S^1\times(-\infty,0],\\
 &\Sigma_2\approx [0,1]\times(-\infty,0],\\
 &\Sigma_3\approx [0,1]\times[0,+\infty),\\
 &\Sigma_0:=\Sigma\setminus(\Sigma_1\cup \Sigma_2\cup \Sigma_3)
 \end{aligned}$$
 and by $u_j:=u|_{\Sigma_j}$.

\begin{center}\includegraphics[width=6cm,height=2cm]{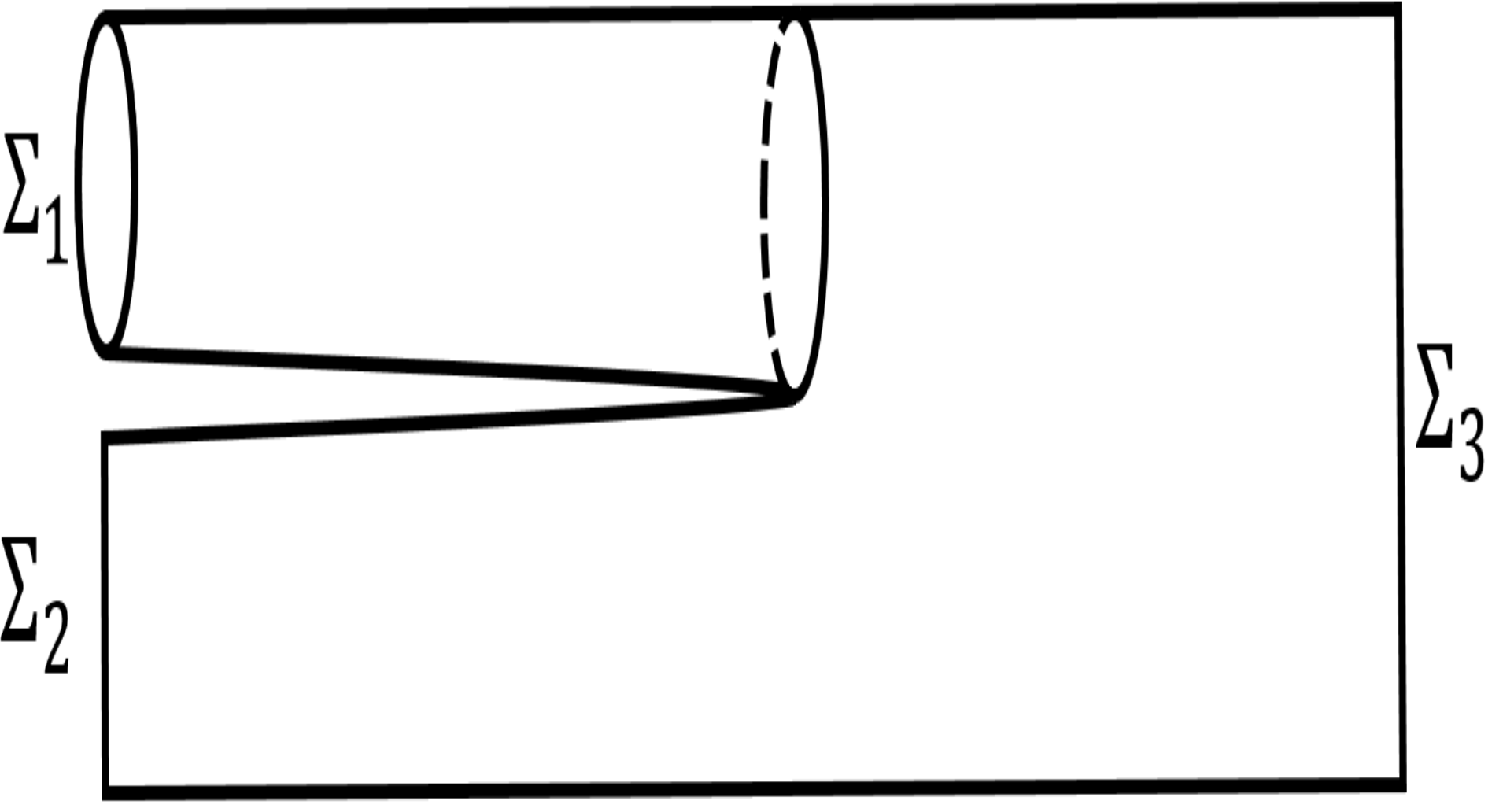}\\
\end{center}
\vspace{3mm}
\begin{center}{\bf Figure 3:}\, Riemannian surface $\Sigma$
\end{center}

\medskip

 Let $\rho_j:\R\to[0,1]$ denote the smooth cut-off functions such that
 $$\rho_1(s)=\rho_2(s)=\begin{cases}1,&s\le -2,\\0,&s\ge -1\end{cases}\quad\rho_3(s):=\rho_1(-s).$$
 For $a\in CF_*(H_1)$, $x\in CF_*(H_2:L)$ and $y\in CF_*(H_3:L)$, denote by
 $\mathcal M(a,x;y)$ the set of all $u$ that satisfy
 $$\left\{\begin{array}{l}
 u:\Sigma\to P,\\
 \pa_su_j+J(\pa_tu_j-X_{\rho_jH_j}\circ u_j)=0,\,j=1,2,3,\\
 \pa_su+J\pa_tu=0,\,\mbox{on}\;\Sigma_0,\\
 u(s,-1),u(s,0^-)\in L,\,s\le 0,\\
 u(s,-1),u(s,1)\in L,\,s\ge 0,\\
 u_1(-\infty,t)=a(t),\\
 u_2(-\infty,t)=x(t),\\
 u_3(+\infty,t)=y(t).
 \end{array}\right.$$
 For generic choices, the set $\mathcal M(a,x;y)$ is a smooth manifold of dimension $\mu_{H_1}^{CZ}(a)+\mu_{H_2}(x)-\mu_{H_3}(y)+n$,
 where $\mu^{CZ}$ denotes the Conley--Zehnder index, and $\mu_{H_j}$ the (corresponding) Maslov index.

\begin{center}\includegraphics[width=7cm,height=3cm]{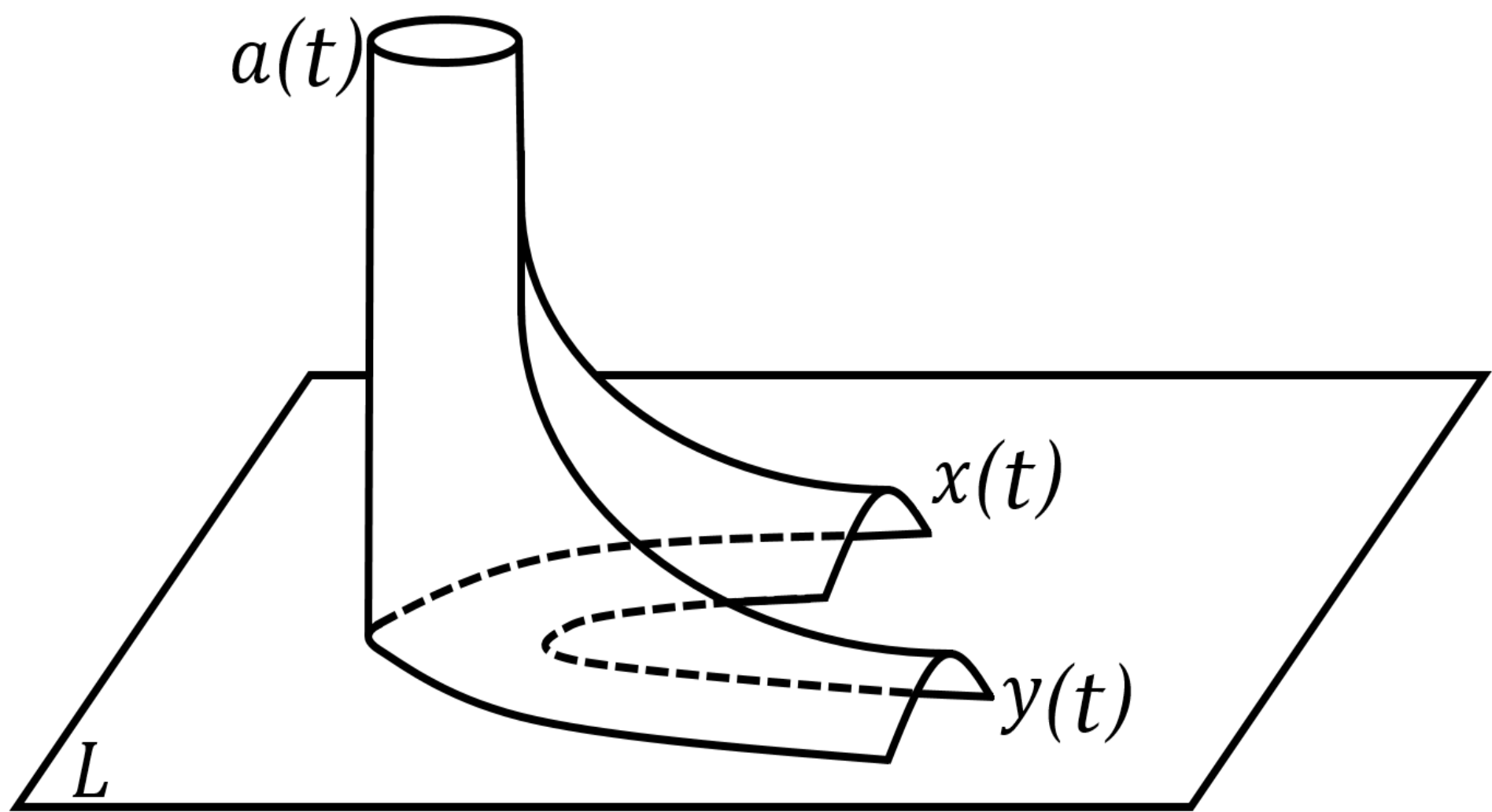}\\
\end{center}
\vspace{3mm}
\begin{center}{\bf Figure 4:}\, Moduli space $\mathcal M(a,x;y)$
\end{center}

\medskip

 Let $n(a,x;y)$ denote the number (modulo $2$) of $\mathcal M(a,x;y)$ in dimension zero. Then we define the map product
$$HF_*(H_1)\otimes HF_*(H_2:L)\to HF_*(H_3:L)$$
by
$$a\circ x:=\sum_yn(a,x;y)y$$ on generators and extend it by (bi)linearity on $HF_*(H_1)\otimes HF_*(H_2:L)$. Using standard cobordism arguments, one can show that $\circ$ descends to the homology level and, when $H_2=H_3$, it defines the product that makes $HF_*(H_2:L)$ a module over $HF_*(H_1)$.

\subsection{Proof of Theorem~\ref{thm:product}}
In order to prove the inequality~(\ref{eq:ineq_circ}) and the Theorem~\ref{thm:product}, we consider, as in~\cite{Oh3} and~\cite{Sc2}, the bundle
$\widetilde{P}\to\Sigma$ whose fiber is isomorphic to $(P,\omega)$ and we fix the trivializations
$$\varphi_j:\widetilde{P}_j:=\widetilde{P}|_{\Sigma_j}\to\Sigma_j\times P$$ for $j=1,2,3$. On each $\widetilde{P}_j$ let
$$\tilde{\omega}_j:=\varphi_j^*(\omega+d(\rho_jH_jdt)).$$ We will use the following theorem by Entov:

\begin{thrm}\label{thm:Entov}\cite{E} There exist a closed two form $\tilde{\omega}$ such that
\begin{itemize}
\item[(1)] $\tilde{\omega}|_{\Sigma_j}=\tilde{\omega}_j$;
\item[(2)] $\tilde{\omega}$ restricts to $\omega$ at each fiber;
\item[(3)] $\tilde{\omega}^{\wedge(n+1)}=0$.
\end{itemize}
\end{thrm}
 Let $\tilde{\omega}$ be as in Theorem~\ref{thm:Entov} and let
$$\Omega_\lambda:=\tilde{\omega}+\lambda\omega_\Sigma,$$ where $\omega_{\Sigma}$ is an area form on $\Sigma$ such that
$\int_{\Sigma}\omega_\Sigma=1$. Choose an almost complex structure $\widetilde{J}$ on $\widetilde{P}$ such that
\begin{itemize}
\item[(1)] $\widetilde{J}$ is $\tilde{\omega}$ compatible on each fiber, so it preserves the vertical tangent space;
\item[(2)] the projection $\pi:\widetilde{P}\to\Sigma$ is $\widetilde{J}-i$ pseudoholomorphic, i.e. $d\pi\circ\widetilde{J}=i\circ d\pi$, where $i$ is a complex structure on $\Sigma$;
\item[(3)] $(\varphi_j)_*\widetilde{J}=i\oplus J_j$, where $J_j(s,t,x):=(\phi_{\rho_jH_j}^t)^*J$.
\end{itemize}
With such a choice, we get that the $\widetilde{J}-$holomorphic section $\tilde u$ over $\Sigma_j$ (or some shorter cylindrical ends, i.e., diffeomorphic to $(-\infty,K_j]\times S^1$, etc.) are precisely the solutions of
\begin{equation}\label{eq:ends}
\pa_su+J\left(\pa u_s-X_{(\rho_jH_j)}\circ u\right)=0.\end{equation}
As in~\cite{Sc2} or~\cite{E} we obtain, for $a\in CF_*(H_1)$, $x\in CF_*(H_2:L)$ and $y\in CF_*(H_1\sharp H_2:L)$
\begin{equation}\label{eq:action_sum}
\int \tilde u^*\tilde\omega=\mathcal{A}_{H_1}(a)+a_{H_2}(x)-a_{H_1\sharp H_2}(y),
\end{equation}
whenever there exists a $\widetilde{J}-$holomorphic section $\tilde u:\Sigma\to\widetilde P$
that satisfies~(\ref{eq:ends}) on fibers. Since $\widetilde{J}$ is $\Omega_\lambda-$ compatible, it holds
$$0\le\int\tilde u^*\Omega_\lambda=\int\tilde u^*\omega+\lambda\int\tilde u^*\omega_\Sigma
=\int\tilde u^*\omega+\lambda\int\omega_\Sigma=\int\tilde u^*\omega+\lambda.$$ Now we use the Entov's result again, that enables us to choose, for any $\delta>0$, a closed two form $\tilde\omega$ such that $\Omega_\lambda$ is symplectic for all $\lambda\ge\delta$ (see~\cite{E} Theorems 3.6.1 and 3.7.4).

Let $\delta>0$, $a\in HF_*(H_1)$, $b\in HF_*(H_2:L)$. Let $\tilde a$ and $x$ be representatives of the class $a$ and $b$ respectively, such that
$$\mathcal A_{H_1}(\tilde a)\le \rho(\PSS^{-1}(a),H_1)+\delta,\quad a_{H_2}(x)\le\sigma(\PSS^{-1}(b),H_2)+\delta.$$
For any $y\in a\circ b$, there exists $u\in\mathcal M(\tilde{a},x;y)$, so we have
$$\begin{aligned}a_{H_1\sharp H_2}(y)&\le \mathcal A_{H_1}(\tilde a)+a_{H_2}(x)+\delta\\
&\le\rho(\PSS^{-1}(a),H_1)+\delta+\sigma(\PSS^{-1}(b),H_2)+\delta+\delta\\
&=\rho(\PSS^{-1}(\alpha),H_1)+\sigma((\PSS^{-1}(b),H_2)+3\delta.
\end{aligned}$$
Since the above inequality is true for all $\delta>0$ and $y$, we conclude
$$\sigma((\PSS^{-1}(a\circ b),H_1\sharp H_2)\le\rho(\PSS^{-1}(a),H_1)+\sigma((\PSS^{-1}(b),H_2),$$
so the Theorem~\ref{thm:product} follows.

\begin{rem} \rm For a smooth submanifold $L$ of $P$ and three Morse functions
$$f_1:P\to\R, \quad f_2, f_3:L\to\R$$ one can define a Morse homology product
$$\cdot:HM_*(P,f_1)\otimes HM_*(L,f_2)\to HM_*(L,f_3)$$ as follows. Let $p_j$ be critical points of $f_j$, for $j=1,2,3$. The set $\mathcal M(p_1,p_2;p_3)$ is defined as the set of all trees $\gamma:=(\gamma_1,\gamma_2,\gamma_3)$ such that
$$\left\{\begin{array}{l}\gamma_1:(-\infty,0]\to P,\; \gamma_2:(-\infty,0]\to L,\; \gamma_3:[0,+\infty)\to L,\\
\dot\gamma_j=-\nabla f_j(\gamma_j),\,j=1,2,3,\\
\gamma_1(-\infty)=p_1,\;\gamma_2(-\infty)=p_2,\;\gamma_3(+\infty)=p_3,\\
\gamma_1(0)=\gamma_2(0)=\gamma_3(0).
\end{array}\right.$$

\begin{center}\includegraphics[width=6cm,height=2.5cm]{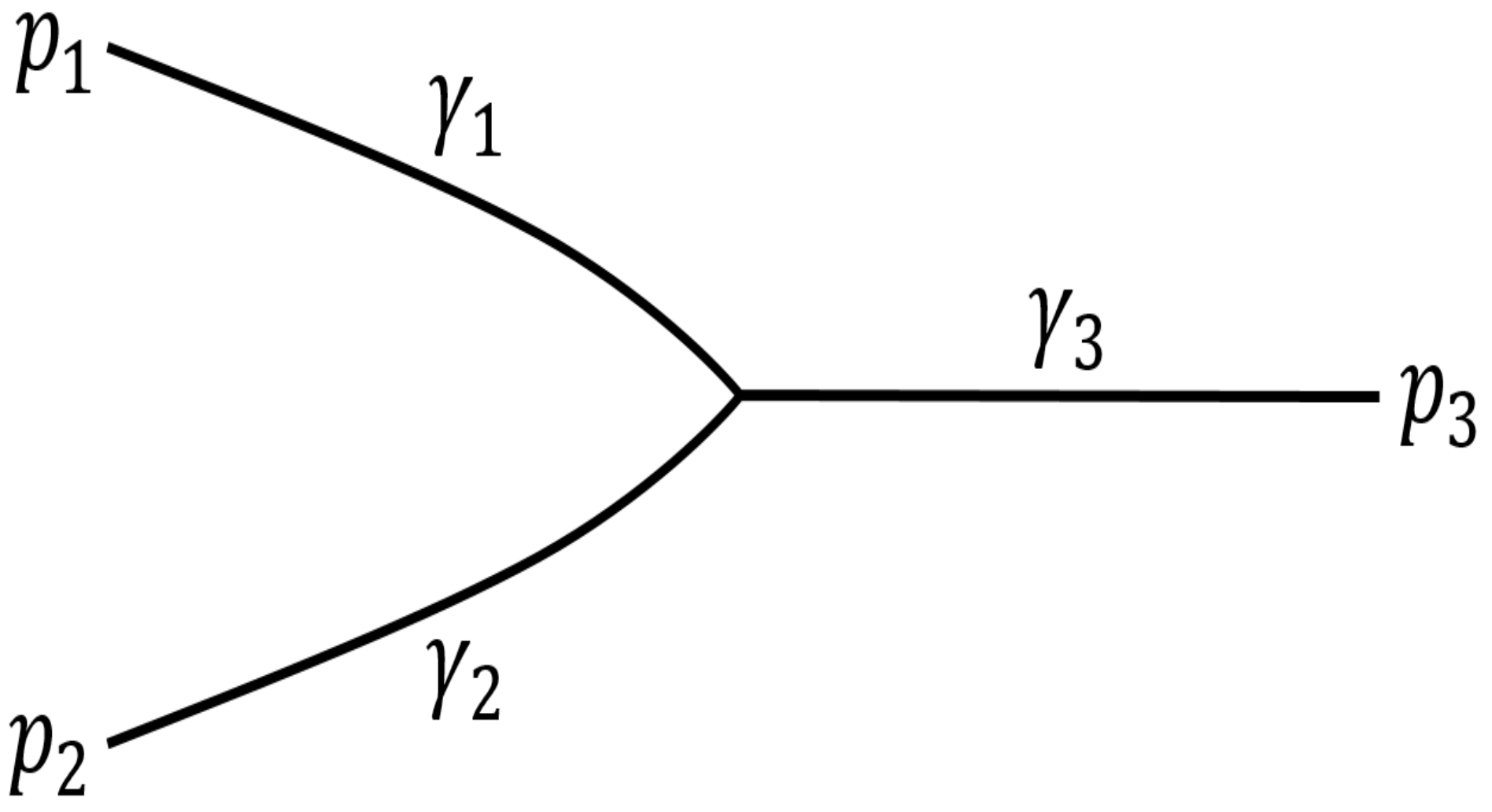}\\
\end{center}
\vspace{3mm}
\begin{center}{\bf Figure 5:}\, The set of trees $\mathcal M(p_1,p_2;p_3)$
\end{center}

\medskip
For generic choices the set $\mathcal M(p_1,p_2;p_3)$ is a smooth manifold of dimension
$$m_{f_1}(p_1)+m_{f_2}(p_2)-m_{f_3}(p_3)-\dim P$$ where $m_f$ is the corresponding Morse index. If $n(p_1,p_2;p_3)$ denotes the number of zero-dimensional component, then the product $\cdot$ is defined as:
$$p_1\cdot p_2=\sum_{p_3} n(p_1,p_2;p_3)p_3$$
on generators.

Using the standard cobordism arguments, one can prove that the product $\circ$ and $\cdot$ commute with PSS type isomorphisms, more precisely,
for $\alpha\in HM_*(P)$, $\beta\in HM_*(L)$, it holds:
\begin{equation}\label{eq:aux1}
\PSS(\alpha\cdot\beta)=\PSS(\alpha)\circ\PSS(\beta),
\end{equation} where $\PSS$ denote both types are PSS-type isomorphisms. Let 
$$a:=\PSS(\alpha),\quad b:=\PSS(\beta).$$
From~(\ref{eq:ineq_circ}) and~(\ref{eq:aux1}) we get
$$\begin{aligned}
&\sigma(\alpha\cdot\beta,H_1\sharp H_2)=\sigma(\PSS^{-1}(a\circ b),H_1\sharp H_2)\le\\
&\rho(\PSS^{-1}(a),H_1)+\sigma(\PSS^{-1}(b),H_2)=
\rho(\alpha,H_1)+\sigma(\beta,H_2),\end{aligned}$$
so
$$\sigma(\alpha\cdot\beta,H_1\sharp H_2)\le \rho(\alpha,H_1)+\sigma(\beta,H_2).$$
\end{rem}


\end{document}